\documentclass[11pt]{article}
\usepackage{amsmath}
\usepackage{amsfonts}
\usepackage{amsthm}
\usepackage{verbatim}

\begin{document}

\def\N{\mathbb{N}}
\def\F{\mathbb{F}}
\def\Z{\mathbb{Z}}
\def\R{\mathbb{R}}
\def\Q{\mathbb{Q}}
\def\H{\mathcal{H}}

\parindent= 3.em \parskip=5pt

\centerline{\bf{A note on hyperquadratic elements}} 
\centerline{\bf{of low algebraic degree}}

\vskip 0.5 cm
\centerline{\bf{by A. Lasjaunias}}
\vskip 0.5 cm
\noindent {\bf{Abstract.}} In different areas of discrete mathematics, a certain type of polynomials, having coefficients in a field $K$ of finite characteristic, has been considered. The form and the degree of these polynomials, here called projective, are simply linked to the 
      characteristic $p$ of $K$. Roots of these projective polynomials are particular algebraic elements over $K$, called hyperquadratic. 
      For a general algebraic element of degree $d$ over $K$, we discuss the possibility of being hyperquadratic. Using a method of differential algebra, we obtain, for particular fields $K=\F_p$, projective polynomials only having polynomial factors of degree 1 or 2.
\vskip 1 cm
\par Let $K$ be a field of positive characteristic $p$ and $r=p^t$ where $t\geq 0$ is an integer. To a given quadruple $(u,v,w,z)$ in $K^4$, such that $uz-vw \neq 0$, we associate a polynomial $H_{K;r}$ (or simply $H$) in $K[X]$, defined by :

    $$H(x)=ux^{r+1}+vx^r+wx+z.$$

   \par These polynomials have been considered long ago, probably first by Carlitz \cite{C}, and studied more recently 
   from an algebraic point of view in a general context by several authors \cite{A},\cite{B}. Following Abyankhar, we call $H$ a projective polynomial over $K$. To be more precise, we can say a projective polynomial of order $t$. We observe that $H(x)=0$ is equivalent to $x=f(x^r)$ where $f$ is a linear fractional transformation defined by $f(x)=(-vx-z)/(ux+w)$. The condition  $uz-vw \neq 0$ makes sure that this linear fractional transformation $f$ is non-trivial and invertible.

     \par If $\alpha \notin K$ is such that there exists a projective polynomial $H$ and we have $H(\alpha)=0$, we say that $\alpha$ is hyperquadratic over $K$. Hence, a hyperquadratic element is a fixed-point of the composition of a linear fractional transformation and of the Frobenius ismorphism $x \rightarrow x^r$. To be more precise, we say that an irrational root (i.e. $\notin K$) of $H_{K,r}$ is an hyperquadratic element of order $t$. Note that a hyperquadratic element over $K$ is a particular algebraic element over $K$ of degree $d$ with $2\leq d \leq r+1$.
   
      \par Since $x \rightarrow x^r$ is an isomorphism in $K$, we have the following: 
      if $\alpha$ is hyperquadratic of order $t$ then we have $\alpha=f(\alpha^r)=f((f(\alpha^r))^r)=g(\alpha^{r^2})$, where $g$
      is another invertible linear fractional transformation. Hence $\alpha$ is also hyperquadratic of order $2t$, and by iteration
      of order $mt$ for all integers $m\geq 1$.  
      \par If $r=1$ (i.e. $t=0$), then $H$ is a polynomial of degree $2$. Hence quadratic elements over $K$ are hyperquadratic elements of order $0$. If $r>1$ and $\alpha$ is algebraic over $K$ of degree $2\leq d \leq 3$, then the four elements $1,\alpha,\alpha^r$ and $\alpha^{r+1}$ in $K(\alpha)$ are linked over $K$. Consequently there exists a polynomial $H$ such that $H(\alpha)=0$ and therefore $\alpha$ is hyperquadratic of any order $t\geq 1$. Accordingly, to be more precise, we define the absolute order of a hyperquadratic element $\alpha$ as the least integer $t$ such there
      is $H$ with $H(\alpha)=0$ and $r=p^t$. Hence a quadratic element over $K$ is hyperquadratic of absolute order $0$ (but also of any order $t\geq 0$) and a cubic element over $K$ is hyperquadratic of absolute order $1$ (but also of any order $t\geq 1$).
      
      \par In this note, for the field $K$, we will only be considering the following two cases. The first case is $K$ finite and consequently $K=\F_q$ where $q$ is a power of a prime $p$. The second case is $K$ being a transcendental extension of a finite field, that is $K=\F_q(T)$ where $T$ is a formal indeterminate. Note that the first case can just be seen as a particular case of the second one. For $K=\F_q$, the sudy of $H$ appears in different works, some more general and others oriented to coding theory (see \cite{BGM},\cite{ST}, \cite{MR1},\cite{MR2},\cite{MR3}). The importance of $H$ in the second case appears in diophantine approximation and continued fractions in function fields over a finite field. The first consideration in this setting, with $K=\F_2(T)$, is due to Baum and Sweet \cite{BS}. For a survey and more references in this area the reader may consult \cite{L2}. As we will see below the study of $H$ in this second case, allows to use methods which bring results also in the first case.
      
      \par Let us consider the case $K=\F_q(T)$. A method to study rational approximation of roots of $H$ in power series fields, based on arguments of differential algebra, was developed. See \cite{L2}, for more precisions and references (note that hyperquadratic elements were first called algebraic of class I). For a short presentation of the arguments developed below, the reader may also consult \cite[p.~260-262]{BL}.  We use formal differentiation in $K$. If $x\in K$ (or a field extension of $K$), we denote by $x'$ the formal derivative of $x$ respect to $T$. 
      If $\alpha$ is algebraic of degree $d$, there is a polynomial $P\in K[x]$  of degree $d$ such that we have $P(\alpha)=0$. By differentiation, we get $\alpha'P'_X(\alpha)+P'_T(\alpha)=0$ and consequently
         $\alpha' \in \F_q(T,\alpha)$. Therefore we get $\alpha'=Q(\alpha)$ where $Q$ is a polynomial of degree less or equal to $d-1$, with coefficients in $\F_q(T)$.
         \par Just to illustrate the above argument, let us consider the simple case $d=2$ : $x$ satisfies $x^2+ax+b=0$ where $a,b \in K$ with $p>2$. 
         Then setting $\Delta=a^2-4b$, through a basic computation the reader may check that we get $\Delta x'=(aa'-2b')x+2ba'-ab'$. 
         We report here below the computation by means of electronic media, applying PARI/GP (This computation can be performed online at https://pari.math.u-bordeaux.fr/gp.html). Given a polynomial $P$, the polynomial $\Delta Q$ is returned (where $\Delta$ is the discriminant of $P$) . Here the derivatives of $a$ and $b$ are denoted by $ap$ and $bp$ respectively.
         \begin{verbatim}       
          ? P=Pol([1,a,b]);Pt=Pol([ap,bp]);
          [U,V,R]=polresultantext(P,P');Q=V*Pt%P
            %1 = (ap*a - 2*bp)*x + (-bp*a + 2*ap*b)
        \end{verbatim}

         \par Returning to the general case, if $\alpha$ is a hyperquadraic element, since $\alpha=f(\alpha^r)$, we get $\alpha'=Q(\alpha)$ with $\deg(Q)\leq 2$ (see \cite[p. 262, Proposition 2.2]{BL}).   
         
         Hence a hyperquadratic element satisfies a Riccati differential equation, in other words we say that it is a differential-quadratic element. Incidently, this shows that, for a general algebraic element of large degree $d$ over $K$, the possibility of being hyperquadratic is remote.

         \par Indeed, from $d\geq 4$ on, the situation is more complex : a general algebraic element of degree $d$ over K may not be differential-quadratic and therefore
         it cannot be hyperquadratic. Starting from this observation, we could ask the following: Given a general polynomial $P$ of degree $d=4$, is there a simple condition on its coefficients such that the root of $P$ is differential-quadratic ?
         The polynomial in its general form, after a translation on $x$, for a characteristic $p>3$, can be written as $P(x)=x^4+ax^2+bx+c$. It was proved that $a^2+12c=0$ is a condition which implies 
         that an eventual root of $P$ is differential-quadratic (see \cite[p. 262]{BL}).
         This can be checked using computer calculations. We write here below the code using PARI/GP as above. The polynomial returned has degree 3 (here as above $ap$, $bp$ and $cp$ stand for the derivatives $a'$, $b'$ and $c'$). 
        \begin{verbatim}
         ?  P=Pol([1,0,a,b,c]);Pt=Pol([ap,bp,cp]);
         [U,V,R]=polresultantext(P,P');Q=V*Pt%P
          
         %1 = (-8*cp*a^3+(4*bp*b+16*ap*c)*a^2+(-6*ap*b^2+32*cp*c)*a+
         (-36*cp*b^2+48*bp*c*b-64*ap*c^2))*x^3+((4*cp*b+16*bp*c)*a^2+
          (-6*bp*b^2-32*ap*c*b)*a+(9*ap*b^3+48*cp*c*b-64*bp*c^2))*x^2
         +(-8*cp*a^4+(4*bp*b+8*ap*c)*a^3+(-4*ap*b^2+48*cp*c)*a^2+
         (-42*cp*b^2+16*bp*c*b-32*ap*c^2)*a+(9*bp*b^3-12*ap*c*b^2-64*cp*c^2))*x
         +((-4*cp*b+8*bp*c)*a^3-4*ap*c*b*a^2+(48*cp*c*b-32*bp*c^2)*a+
         (-27*cp*b^3+36*bp*c*b^2-48*ap*c^2*b)).
      \end{verbatim}
      
      And finally, after the substitution $c=-a^2/12$ and the one obtained by differentiation, we observe that the leading coefficient of $Q$ vanishes.
      
      \begin{verbatim}
        ? substvec(Q,[c,cp],[-a^2/12,-a*ap/6])
       %2 = (-16/9*bp*a^4+8/3*ap*b*a^3-6*bp*b^2*a+9*ap*b^3)*x^2
          +(32/27*ap*a^5+8/3*bp*b*a^3+4*ap*b^2*a^2+9*bp*b^3)*x
          +(-8/9*bp*a^5+4/3*ap*b*a^4-3*bp*b^2*a^2+9/2*ap*b^3*a).
      \end{verbatim}
      
         \par Then a natural question arises: under the condition $a^2+12c=0$, may a solution of $P$ be hyperquadratic ? 
         The answer is positive. Indeed, in \cite[p.~35-38]{L1} with a limitation on the size of the prime $p$, and in \cite{BL} without limitation, the following was proved:
         For $p>3$ and $p \equiv i \mod 3 (i=1, 2)$, $a,b \in K$, the polynomial $P(x)=x^4+ax^2+bx-a^2/12$ divides a projective polynomial of order $i$. Just to briefly illustrate this: if $p=7$ and $a,b \in K$, we have
         $$ax^8+3bx^7+4b(b^2+4a^3)x+2a^2(b^2+a^3)=$$
         $$(x^4+ax^2+bx+4a^2)(ax^4+3bx^3+6a^2x^2+3abx+4(b^2+a^3)).$$  
           
        \par The existence of such a simple condition, on the coefficients of the polynomial $P$, implying it to divide a projective polynomial remains
        somehow mysterious. Thus, we decided to investigate the case $d=5$, searching for eventual differential-quadratic elements.
        After a translation on $x$, the general form of $P$ would be $P=x^5+ax^3+bx^2+cx+d$ for $p>5$. The polynomial $Q$, such that $x'=Q(x)$, would be of degree 4:      $Q=b_4x^4+b_3x^3+b_2x^2+b_1x+b_0$. Hence we need to check the coefficients $b_4$ and $b_3$, trying to find which conditions
            on $a,b,c$ and $d$ would make them both vanish. The computations to obtain the 5 coefficients of $Q$, have been performed as above using PARI. However, the situation appears too intricated due to the number 4 of coefficients in $P$. To simplify, we decided to check the simpler case of $P$ having no term of degree 3. Our goal was to obtain a hyperquadratic element algebraic of degree 5. However, we were unsuccessful. We could only obtain very partial results, bringing more questions than answers, which we expose here below.
           
            \par We consider $P=x^5+ax^2+bx+c$ with $a,b$ and $c$ in $\F_q(T)$  and $p>5$. After a thoroughful examination of the coefficients $b_4$ and $b_3$, we observed the following. Under a couple of particular sufficient conditions $(C_1)$ and $(C_2)$ on the three coefficients $a,b$ and $c$, we have $b_3=0$ and $b_4=0$.    
            These conditions are the following:     
            $$(C_1)\quad  18a^3+325bc=0 \quad \text{and}\quad   (C_2)\quad  5b'c=4c'b.$$ 
                   
            \par Hence, if $(C1)$ and $(C2)$ are satisfied then a solution of $P$ is differential-quadratic. (We checked the other coefficients $b_2,b_1$ and $b_0$ and we observed that we also have $b_2=b_0=0$ !). 
            The question is: under conditions $(C_1)$ and $(C_2)$ could this solution be hyperquadratic ?
            We could only give a very partial answer to this question. Note that condition $(C_2)$ can be written as $(b^5/c^4)'=0$ if $c \neq 0$. 
            We introduce the condition $(C_3)\quad b^5=2c^4$. Note that $(C_3)$ implies $(C_2)$.
            Our result is the following: Let $K=\F_p$, $P$ as above and $a,b$ and $c$ satisfying $(C_1)$ and $(C_3)$. Then, if $p=11$ or $p=17$,  $P$ divides a projective polynomial $H$ of order 1.
            \par  This was obtained by direct computations. Amazingly, the attempt to obtain the same for other prime numbers was unsuccessful. Moreover, in all these cases, the polynomial $P$ is splited in the same form $2^2*1$ (two factors of degree 2 and one of degre 1), while the corresponding polynomial $H$ has $(p+3)/2$ factors and it is is splited in the form $2^{(p-1)/2}*1*1$.
                        
            \par First we show how the three coefficients of $P$ have been obtained satisfying the above conditions. Once $P$ is chosen, to possibly obtain  the polynomial $H$, it is enough to compare the remainders modulo $P$ of $x^{r+1}$ and $x^r$ respectively and then to check whether a linear combination of these ones forms a polynomial of degree 1. 
            
            \par If $p=6k+5$, we observe that the map $x\rightarrow x^3$ is one to one in $\F_p$. We denote the inverse map by $x\rightarrow cr(x)$ and we simply have  $cr(x)=x^{-2k-1}$ in $\F_p^*$.
            For $p\neq 5,13$, we set $u=2(18/325)^4 \in \F_p^*$. Let us consider the triple $(a,b,c) \in \F_p^3$ where 
            $$p=11,17 \quad a\in \F_p^*\quad b=cr(a^4cr(u))\quad \text{and} \quad c=(-18a^3)/(325b).$$ It is easy to check that the triple $(a,b,c)$ satisfies conditions $(C_1)$ and $(C_3)$. Each triple $(a,b,c)$ will correspond to a polynomial $P$, hence we have $10+16$ possible cases.

            \par Here below, in two tables corresponding to the cases $p=11$ and $p=17$ respectively, we describe the polynomials $P$ and $H$ in $\F_p[X]$ such that $P$ divides $H$. In these tables the polynomials $P=x^5+ax^2+bx+c$ and $H=ux^{p+1}+vx^p+wx+z$ , where $a,b,c,u,v,w$ and $z \in \F_p$, are respectively represented by the tuples $(a,b,c)$ and $(u,v,w,z)$. Moreover $H$ is defined up to a constant factor and consequently we may choose it to be unitary.

\vskip 0.5 cm
\centerline{Table 1: $p=11$}
\begin{center}
   \begin{tabular}{| c | c | c | c | c | }
     \hline
      $P$ & $H$ && $P$ & $H$ \\ \hline
     (1,7,9) & (1,7,7,2)&& (6,6,2) & (1,1,1,7) \\ 
     (2,10,2) &  (1,5,5,10)&& (7,8,2) & (1,9,9,6) \\
     (3,2,9) & (1,8,8,8) &&(8,2,2) & (1,3,3,8)\\
      (4,8,9) & (1,2,2,6)&&(9,10,9) & (1,6,6,10)\\
      (5,6,9) & (1,10,10,7)&&(10,7,2) & (1,4,4,2) \\
     \hline
        
   \end{tabular}
   \end{center}

   \vskip 0.5 cm
\centerline{Table 2: $p=17$}
\begin{center}
   \begin{tabular}{| c | c | c | c | c | }
     \hline
      $P$ & $H$ && $P$ & $H$ \\ \hline
     (1,15,13) & (1,13,13,3)&& (9,2,9) & (1,8,8,12) \\
     (2,2,15) & (1,2,2,5)&& (10,8,14) & (1,5,5,10) \\
     (3,9,7) & (1,6,6,11)&& (11,8,5) & (1,3,3,7) \\
     (4,15,16) & (1,16,16,14)&& (12,9,6) & (1,10,10,6) \\
     (5,9,11) & (1,7,7,6)&& (13,15,1) & (1,1,1,14) \\
     (6,8,12) & (1,14,14,7)&& (14,9,10) & (1,11,11,11) \\
     (7,8,3) & (1,12,12,10)&& (15,2,2) & (1,15,15,5) \\
     (8,2,8) & (1,9,9,12)&& (16,15,4) & (1,4,4,3) \\
      \hline
        
   \end{tabular}
   \end{center}
\vskip 0.5 cm

\noindent {\bf{Aknowledgements. }} We would like to thank Bill Allombert for his skillful advices on computer programming and his help in using  PARI system.

\vskip 0.5 cm
\begin{tabular}{ll}Alain LASJAUNIAS
\\Rue du Livran \\L\'eognan 33850, France \\E-mail: lasjauniasalain@gmail.com\\\end{tabular}

\end{document}